\documentclass[12pt]{article}

\usepackage{amsfonts,amssymb,epsfig}
\usepackage{amsmath}
\usepackage{amsthm}

\newtheorem{thm}{Theorem}
\newtheorem*{thm*}{Theorem}
\newtheorem{lm}{Lemma}

\def\pf{ \noindent {\bf Proof: \  }}
\def\endpf{\qed \medskip} 
\def\e{\varepsilon}

\begin{document}

\title{Stabilizing isomorphisms from $\ell_{p}\left(\ell_{2}\right)$ into $L_{p}\left[0,1\right]$}

\author{Ran Levy\thanks{This is part of the first author MSc thesis written at the Weizmann Institute under the supervision of the second author}\ \  and Gideon Schechtman\thanks{supported in part by the Israel Science Foundation and by the Israel-US binational Science Foundation}}

\maketitle

\abstract

Let $1<p\not=2<\infty$, $\e>0$ and let $T:\ell_p(\ell_2)\overset{into}{\rightarrow}L_p[0,1]$ be an isomorphism
Then there is a subspace $Y\subset \ell_p(\ell_2)$
$(1+\e)$-isomorphic to $\ell_p(\ell_2)$ such that: $T_{|Y}$
is an $(1+\e)$-isomorphism and $T\left(Y\right)$ is $K_p$-complemented
in $L_{p}\left[0,1\right]$, with $K_p$ depending only on $p$. Moreover, $K_p\le (1+\e)\gamma_p$ if $p>2$ and $K_p\le (1+\e)\gamma_{p/(p-1)}$
if $1<p<2$, where $\gamma_r$ is the $L_r$ norm of a standard Gaussian variable.

\section{Introduction}

Let $B$ be one of the 5 ``classical" subspaces of $L_p=L_p[0,1]$; by these we mean $\ell_p$, $\ell_2$, $\ell_p\oplus_p\ell_2$, $\ell_p(\ell_2)$ and $L_p$ itself. Here $\ell_p(\ell_2)$ is the space of (say, real) matrices $A=(a_{i,j})_{i,j=1}^\infty$ with norm $\|A\|_{\ell_p(\ell_2)}=(\sum_{j=1}^\infty(\sum_{i=1}^\infty a_{i,j}^2)^{p/2})^{1/p}$.
It is well known that these five spaces isometrically embed in $L_p$, $1\le p<\infty$, and, for $1<p<\infty$, have embeddings which are complemented. (For $p=1$ this last statement holds only for $L_1$ and $\ell_1$.)

It was known for some time that if $X$ is any subspace of $L_p$, $1<p<\infty$, isomorphic to one of these $B$ then there is a subspace of $X$ isomorphic to $B$ and complemented in $L_p$. For $\ell_p$ see \cite{KP}, for $\ell_2$, \cite{PR}.  The case $\ell_p\oplus_p\ell_2$ follows from these two results. The quite complicated case of $L_p$ was proved in \cite{JMST} (and for $p=1$ perviously in \cite{ES}). The case of $\ell_p(\ell_2)$ can be proved using a variation of the method of \cite{JMST} (and there is also a much simpler proof for $p>2$) and was known to the second named author for a long time but not published (the simpler proof for $p>2$ is included in \cite{HOS}).

Recently, there were three paper which address this property for some of the spaces $B$ above and related questions again. This was done mostly because some strengthening of this property was needed for other purposes. Firstly, Haydon, Odell and Schlumprecht proved in \cite[Theorem 6.8]{HOS} that, for $p>2$, any subspace of $L_p$ isomorphic to $\ell_p(\ell_2)$ contains a subspace $(1+\e)$-isomorphic to $\ell_p(\ell_2)$ and complemented in $L_p$ by means of a projection of norm at most $(1+\e)\gamma_p$, $\gamma_p$ being the $L_p$ norm of a standard Gaussian variable. Their proof uses the fact that a similar theorem holds for the space $\ell_2$ in $L_p$; i.e., any subspace of $L_p$, $p>2$, isomorphic to $\ell_2$ contains a subspace $(1+\e)$-isomorphic to $\ell_2$ and complemented in $L_p$ by means of a projection of norm at most $(1+\e)\gamma_p$. This later deep fact, hidden already in Aldous' \cite{Ald}, was recently given a simpler proof by Alspach \cite{Als}.
The third paper is by Dosev, Johnson and the second named author \cite[Theorem 3.4]{DJS} from which the following follows:\\
 {\em For each $1<p<\infty$ there is a constant $K_p$, depending only on $p$ such that if $T:L_p\to L_p$ is an isomorphism (into) then there is a subspace $X$ of $L_p$ $K_p$-isomorphic to $L_p$ such that $T_{|X}$ is a $K_p$-isomorphism and $T(X)$ is $K_p$-complemented in $L_p$.}\\
We remark that a similar theorem for $p=1$ (with the constant $K_1$ arbitrarily close to $1$) is due to Enflo and Starbird \cite{ES}; see also \cite{R} for a somewhat simpler exposition.

The main purpose of the current paper is to prove a similar theorem for $\ell_p(\ell_2)$.
\begin{thm}\label{thm:main}
Let $1<p\not=2<\infty$, $\e>0$ and let $T:\ell_p(\ell_2)\overset{onto}{\rightarrow}X$ be an isomorphism
where $X\subset L_{p}\left[0,1\right]$. Then there is a subspace $Y\subset \ell_p(\ell_2)$
$(1+\e)$-isomorphic to $\ell_p(\ell_2)$ such that $T_{|Y}$
is an $(1+\e)$-isomorphism and $T\left(Y\right)$ is complemented
in $L_{p}\left[0,1\right]$ by means of a projection of norm at most $(1+\e)\gamma_p$ if $p>2$ and $(1+\e)\gamma_{p/(p-1)}$
if $1<p<2$.
\end{thm}

For $p>2$ this theorem follows easily from the result of \cite{HOS} mentioned above so the main innovation here is the case $1<p<2$. However, since the addition needed to present a uniform proof for $1<p<2$ and $p>2$ is minimal, we shall prove both cases. The proof is very much in the spirit of \cite{DJS} but we wrote it in such a manner that one does not need to refer to that paper.

\section{Preliminaries}

The proofs below will assume familiarity with basic techniques of Banach space theory. In particular techniques related to bases. We shall use freely notions like unconditional bases, block bases, small perturbations of bases, gliding hump arguments and similar notions. They can all be found in the first chapter of \cite{LT-I}.

Recall that the Haar system is the following sequence of functions on $[0,1]$:\\
 $h_{0,0}(t)\equiv 1$ and, for
$n = 0, 1,\ldots$ and $i=1,2,\ldots, 2^n$,
\[
h_{n,i}(t) = \left\{ \begin{array}{ll}
1 & \textrm{if}\,\, t\in ((2i-2)2^{-(n+1)}, (2i-1)2^{-(n+1)})\\
-1
& \textrm{if}\,\,t\in ((2i-1)2^{-(n+1)}, 2i2^{-(n+1)})\\
0 & \textrm{otherwise}
\end{array} \right.
\]
This system forms an unconditional basis for $L_p=L_p[0,1]$ for each $1<p<\infty$ (but not in $L_1$ in which it is, in its natural order, a non-unconditional Schauder basis) see e.g. \cite{LT-II}. We denote by $H_p$ its unconditional constant; i.e.,
\[
\|\sum a_{n,i}h_{n,i}\|_p\le H_p\|\sum \e_{n,i}a_{n,i}h_{n,i}\|_p
\]
for any sequence of coefficients $\{a_{n,i}\}$ and any sequence of signs $\{\e_{n,i}\}$. (We shall use the real field although all the arguments easily carry
over to the complex field with minimal changes, one of which should be in this definition). The best constant $H_p$ is known and is of order $\max\{p,p/(p-1)\}$.

Recall Khinchine's inequality; For $1\le p<\infty$,
\[
A_p(\sum_{i=1}^n a_i^2)^{1/2}\leq ({\rm Ave}_{\pm}|\sum_{i=1}^n \pm a_i|^p)^{1/p}\leq B_p
(\sum_{i=1}^n a_i^2)^{1/2}
\]
for all $n$ and all coefficients $\{a_i\}_{i=1}^n$. The best constants $A_p$ and $B_p$ are known and in particular $A_p$ is between $2^{-1/2}$ and $1$ for $1\le p\le 2$ and $1$ for $p>1$ and $B_p$ is $1$ for $1\le p\le 2$ and of order $p^{1/2}$ for $p>2$.

We shall make an intensive use of the square function with respect to the Haar system. For $f\in L_p$, $1<p<\infty$, with expansion $f=\sum a_{n,i}h_{n,i}$ we denote its square function with respect to the Haar system by
\[
S(f)=(\sum a_{n,i}^2h_{n,i}^2)^{1/2}.
\]
The unconditionality of the Haar system and Khinchine's inequality easily imply that
\[
H_p^{-1}A_p\|S(f)\|_p\leq \|f\|_p\leq H_p B_p\|S(f)\|_p,\ \ 1<p<\infty.
\]

Equi-integrability of some sets of functions will play an important role in the sequel. Recall that a set $F$ of Lebesgue integrable functions on $[0,1]$ is said to be {\em equi-integrable} if for all $\e>0$ there is a $\delta>0$ such that for every subset $A$ of $[0,1]$ of measure at most $\delta$, $\int_A|f|dt<\e$ for all $f\in F$.
Equivalently, if For for all $\e>0$ there is a positive $R$ such that $\int|f|{\bf 1}_{|f|>R}dt<\e$ for all $f\in F$.

For $0<r<\infty$ the set $F$ is {\em $r$-equi-integrable} if $\{|f|^r\ ;\ f\in F\}$ is equi-integrable.

Finally, by a {\em $K$-isomorphism} we mean a linear map $T$ from one normed space $X$ into another $Y$ such that $A^{-1}\|x\|\le\|Tx\|\le B\|x\|$ for all $x\in X$ with $AB\le K$. In particular, for small $\e$ a $(1+\e)$-isomorphism $T$ does not necessarily almost preserve the norm of each $x$ (as is sometimes assumed in other places) but of course some multiple of $T$ does.

\section{Stabilazing embeddings of $\ell_{2}$ into $L_p$}\label{sec:3}

In this section we consider the analogue of Theorem \ref{thm:main} for the space $\ell_2$ instead of $\ell_p(\ell_2)$.
As we indicated above this theorem is known although
not simple, especially if one wants to achieve the best constants. (For a somewhat weaker Theorem, in terms of the constants achieved, see Theorem 3.1 in \cite{PR}.) The purpose of this section is to survey its proof and point the reader to the relevant references.

\begin{thm*}\label{thm:l2} Let $\e>0$ and let $T:\ell_2\rightarrow L_{p}\left[0,1\right]$,
$1\leq p<\infty$, be an isomorphism. Then there is an infinite dimensional
subspace $X\subseteq \ell_2$ such that $T_{|X}$ is an $\left(1+\e\right)$-isomorphism
and, for $1<p<\infty$, $TX$ is $\left(1+\e\right)\gamma_{p}$-complemented
in $L_{p}$. Here, for $p>2$, $\gamma_{p}$ is the $L_{p}$ norm
of a standard Gaussian variable and, for $1<p<2$, the $L_{\frac{p}{p-1}}$
norm of such a variable.
\end{thm*}
Note first that in order to prove the first part of the theorem,
the existence of a $\left(1+\e\right)$-isomorphism on $X\subseteq \ell_2$,
it is enough to prove that $TX$ contains a $Y$ which is $\left(1+\e\right)$-isomorphic
to $\ell_2$. Indeed if this is the case let $\left\{ e_{i}\right\} _{i=1}^{\infty}\subseteq X$
be an orthonormal basis of $X$ and $\left\{ f_{i}\right\} _{i=1}^{\infty}\subseteq TX$
a basis $\left(1+\e\right)$-equivalent to an orthonormal basis
of $\ell_2$. Since $\left\{ Te_{i}\right\} _{i=1}^{\infty}$ is a
weakly null sequence in $TX$ for which $\left\{ \left\Vert Te_{i}\right\Vert \right\} _{i=1}^{\infty}$
is bounded and bounded away from zero, we may find a subsequence $\left\{ e_{i_{j}}\right\} _{j=1}^{\infty}$
such that $\left\{ Te_{i_{j}}\right\} _{j=1}^{\infty}$ is as small perturbation
as we wish of a block basis of $\left\{ f_{i}\right\} _{i=1}^{\infty}$
and $\left\{ \left\Vert Te_{i_{j}}\right\Vert \right\} _{j=1}^{\infty}$
is as close as we want to a constant sequence. $X_{0}=span\left\{ e_{i_{j}}\right\} _{j=1}^{\infty}$
is then the subspace we are after. That $TX$ contains a subspace
$\left(1+\e\right)$-isomorphic to $\ell_2$ is by now well known
and follows from the stability of $L_{p}$ (see \cite{KM}).

We are left with the problem of complementation, and especially
the norm of the best projection. For $p>2$ the theorem is specifically
stated and proved in \cite{HOS}. There is a simpler proof in \cite{Als}.
For $1\leq p<2$ this follows, as we shall indicate momentarily from
\cite{Ald}. This is not at all an easy paper to follow and it would
be nice if somebody finds an easier proof maybe a-la-\cite{Als}.
We shall only sketch how to get the result from \cite{Ald} and then give a much simpler argument which however gives a somewhat weaker estimate.

Let $Y\subseteq L_{p}\left[0,1\right]$, $1<p<2$, be isomorphic
to $\ell_2$. We would like to find a subspace $Y_{0}\subseteq Y$
$\left(1+\e\right)$-isomorphic to $\ell_2$ and $\left(1+\e\right)\gamma_{\frac{p}{p-1}}$-complemented
in $L_{p}\left[0,1\right]$.

By \cite{KP} the unit
ball of $Y$ is $p$-equi-integrable and the $L_{1}$ and $L_{p}$
norms are equivalent on $Y$. From here on we shall use the notations
of \cite{Ald}. By the combination of Proposition 3.9 and Theorem
3.10 there, there is a uniformly integrable sequence $X_{n}$ in $Y$ such
that $i\left(X_{n}\right)\overset{wm}{\rightarrow}\sigma\left(q,\alpha\right)$
for some $1<q\leq2$. Recall that $i\left(X_{n}\right)$ is the random
measure $\delta_{X_{n}}$ and that for random measures $\xi_{n}$, $\xi$,
$\xi_{n}\overset{wm}{\rightarrow}\xi$ denotes $\left\langle f,\xi_{n}\right\rangle \overset{w}{\rightarrow}\left\langle f,\xi\right\rangle $
for all $f\in C\left(\mathbb{R}\right)$ (and $\overset{w}{\rightarrow}$
denotes weak convergence in $L_{1}$). Finally, for a function $\alpha\geq0$,
$\sigma\left(q,\alpha\right)$ is the random measure whose characteristic
function is $e^{-\alpha^{q}\left|t\right|^{q}}$; i.e., it is a mixture
of symmetric $q$-stable random variables.

In our case, by Proposition 3.11 of \cite{Ald} the only possible value for $q$
is $q=2$. So we get a sequence $X_{n}$ which tends in some sense
to a mixture of Gaussian variables. Proposition 3.11 and its proof
then say that some subsequence of $X_{n}$ is $\left(1+\e\right)$
equivalent to the unit vector basis of $\ell_2$. The proof really gives more:
some subsequence of $X_{n}$ is arbitrarily close, in $L_{p}$ norm
to a sequence of the form $\alpha Z_{n}$ where, given $\alpha$,
$Z_{n}$ are i.i.d $N\left(0,1\right)$. This means that after a change of density we may assume that
$\{X_{n}\}$ is a small perturbation in the $L_{p}$ norm of a sequence
of i.i.d $N\left(0,1\right)$ variables $Z_{n}$. By a change of density we mean an operator of the form $T_\varphi:L_p\to L_p([0,1],\varphi dt)$, $T_\varphi f=\frac{f}{\varphi^{1/p}}$  where $\varphi$ is a density which is strictly positive on the union of the supports of the $X_n$-s. It is thus enough
to show that the span of such a sequence is $\gamma_{\frac{p}{p-1}}$
complemented in $L_{p}$. Since $\left\{ Z_{n}\right\} _{n=1}^{\infty}\subseteq L_{\frac{p}{p-1}}$
it is clearly enough to show that the orthogonal projection $P:L_{2}\overset{onto}{\rightarrow}\left[Z_{n}\right]$
has norm $\gamma_{\frac{p}{p-1}}$ when considered as an operator
on $L_{\frac{p}{p-1}}$. Since for all $f\in L_{\frac{p}{p-1}}$ (so
also $f\in L_{2}$) $Pf$ is a Gaussian variable $\left\Vert Pf\right\Vert _{\frac{p}{p-1}}=\gamma_{\frac{p}{p-1}}\left\Vert Pf\right\Vert _{2}\leq\gamma_{\frac{p}{p-1}}\left\Vert f\right\Vert _{2}\leq\gamma_{\frac{p}{p-1}}\left\Vert f\right\Vert _{\frac{p}{p-1}}$.\\
This concludes the (admitadly very rough) sketch of the proof of Theorem \ref{thm:l2}

\medskip

We now present a sketch of a {\bf proof of the complementation} part for $1<p<2$
which does not use \cite{Ald}, but gives a somewhat worse constant
than $\left(1+\e\right)\gamma_{\frac{p}{p-1}}$. Let $\left\{ f_{n}\right\} _{n=1}^{\infty}$
be a normalized sequence in $TX$ which is $\left(1+\e\right)$
equivalent to the unit vector basis of $\ell_2$. Let $g_{n}=\left|f_{n}\right|^{p-1}signf_{n}$
so that $\left\Vert g_{n}\right\Vert _{\frac{p}{p-1}}=1$ and $\left\langle g_{n},f_{n}\right\rangle =1$.
Passing to a subsequence we may assume $g_{n}\overset{w}{\rightarrow}g$.
Since $f_{n}\overset{w}{\rightarrow}0$, passing to another subsequence
we may assume that $\left\{ f_{n}\right\} _{n=1}^{\infty}$ and $\left\{ g_{n}-g\right\} _{n=1}^{\infty}$
are arbitrarily close to a biorthogonal sequence; i.e. \[
\sum_{n,m}\left|\left\langle f_{n},g_{m}-g\right\rangle -\delta_{nm}\right|<\e\]
Note that $1-\e\leq\left\Vert g_{n}-g\right\Vert \leq2$ for
all $n$ (the lower bound follows from $\left\langle f_{n},g_{n}-g\right\rangle \geq1-\e$).
We may also assume, by passing to a further subsequence, that $g_{n}-g$
is a martingale difference sequence (or a block basis of the Haar
system). Consequently, for all coefficients $\left\{ a_{n}\right\} _{n=1}^{\infty}$
\[
\begin{array}{rcl}
\left\Vert \sum a_{n}\left(g_{n}-g\right)\right\Vert _{\frac{p}{p-1}} & \leq & H_{p}\mathbb{E}\left\Vert \sum\pm a_{n}\left(g_{n}-g\right)\right\Vert _{\frac{p}{p-1}}\\
 & \leq & 2H_{p}B_{p}\left(\sum a_{n}^{2}\right)^{1/2}\end{array}\]
where $H_{p}$ is the unconditionality constant of the Haar basis
and $B_{p}=\gamma_{\frac{p}{p-1}}$ is the type 2 constant of $L_{\frac{p}{p-1}}$.
Define now $P:L_{p}\rightarrow [f_{n}] $ by $Pf=\sum_{n=1}^{\infty}\left\langle f,g_{n}-g\right\rangle f_{n}$
then it is easily seen that $P$ is a projection of norm $\leq2K_{p}B_{p}\left(1+\e\right)$.

\section{Proof of the main result}

Here we shall prove Theorem \ref{thm:main}. We shall denote by $\left\{ e_{i,j}\right\} _{i,j=1}^{\infty}$ the
canonical basis of $\ell_p(\ell_2)$ i.e.,
\[
\left\Vert \sum_{i,j=1}^{\infty}a_{i,j}e_{i,j}\right\Vert =\left(\sum_{j=1}^{\infty}\left(\sum_{i=1}^{\infty}a_{i,j}^{2}\right)^{p/2}\right)^{1/p}
\]
for all $\left\{ a_{i,j}\right\} _{i,j=1}^{\infty}\subseteq\mathbb{R}$.
By passing to a subsequence in
each column of $\left\{ e_{i,j}\right\} _{i,j=1}^{\infty}$, a gliding
hump argument (applied in the order $\left(11\right),\left(12\right),\left(21\right),\left(13\right),\left(22\right),\left(31\right),...$
of the indices) and a simple perturbation argument, we can assume
that for some infinite subsequences $N_{j}\subseteq N$, $\left\{ Te_{i,j}\right\} _{j=1,i\in N_{j}}^{\infty}$
is a block basis of the Haar system in $L_{p}$. By that we mean that
if the perturbed operator satisfies the conclusion of the theorem so
does the original one. Also, since $\left\{ e_{i,j}\right\} _{j=1,i\in N_{j}}^{\infty}$
spans an isometric $\ell_p(\ell_2)$, and we are anyhow interested
only in a subspace of $\ell_p(\ell_2)$ isometric to $\ell_p(\ell_2)$,
we may also assume that $N_{j}=\mathbb{N}$ for all $j$; i.e., $\left\{ Te_{i,j}\right\} _{i,j=1}^{\infty}$
is a block basis of the Haar system. We shall assume that from now
on.

Given a finite $E\subseteq\mathbb{N}$ and $i\in\mathbb{N}$ set
$v_{i}\left(E\right)=S\left(\sum_{j\in E}Te_{i,j}\right)$. The next
two lemmas are basically taken from \cite{DJS}.
We repeat the proofs for completeness.

\begin{lm}
For all finite $E\subseteq\mathbb{N}$ and $1<p<2$ the convex hull
of $\left\{ v_{i}^{2}\left(E\right)\right\} _{i\in\mathbb{N}}$ is
$p/2$-equi-integrable; i.e., the set of $p/2$ powers of
functions in the convex hull of $\left\{ v_{i}^{2}\left(E\right)\right\} _{i\in\mathbb{N}}$
is equi-integrable.
\end{lm}

\pf
Fix a finite $E\subseteq \mathbb{N}$ and write $v_{i}=v_{i}\left(E\right)$.
Assume that the convex hull of $\left\{ v_{i}^{2}\right\} _{i\in\mathbb{N}}$
is not $p/2$-equi-integrable. Then, there exists $\e_{0}>0$,
a sequence $\{ u_{j}^{2}\} _{j\in\mathbb{N}}$ of disjoint
convex blocks of $\left\{ v_{i}^{2}\right\} _{i\in\mathbb{N}}$ (i.e.
$u_{j}^{2}=\sum_{i\in\sigma_{j}}\alpha_{i,j}^{2}v_{i}^{2}$ where
$\sigma_{j}$ are disjoint subsets of $\mathbb{N}$ and $\sum_{i\in\sigma_{j}}\alpha_{i,j}^{2}\leq1$)
and disjoint subsets $B_{j}$ such that $\int_{B_{j}}|u_{j}|^p>\e_{0}$
for all $j\in\mathbb{N}$. The fact that we can choose the sequence $\{u_j^2\}$ to be disjointly supported with respect to $\{v_i^2\}$ follows from the easy fact that
the convex hull of a finite
subset of $\{ v_{i}^2\} _{i=1}^{\infty}$ is $p/2$-equi-integrable.

Now, For all $\left\{ a_{m}\right\} _{m=1}^{\infty}\in \ell_2$,
we obtain the following inequality, contradicting $p<2$.
\[
\begin{array}{rcl}
\left(\sum_{m=1}^{\infty}a_{m}^{2}\right)^{1/2} & = & \left|E\right|^{-1/p}\left(\sum_{j\in E}\left(\sum_{m=1}^{\infty}a_{m}^{2}\right)^{p/2}\right)^{1/p}\\
 & \geq & \left|E\right|^{-1/p}\left\Vert \sum_{m=1}^{\infty}a_{m}\sum_{i\in\sigma_{m}}\alpha_{i,m}\sum_{j\in E}e_{i,j}\right\Vert _{\ell_p(\ell_2)}\\
 & \geq & \left|E\right|^{-1/p}\left\Vert T\right\Vert ^{-1}H_{p}^{-1}A_{p}\left\Vert S\left(\sum_{m=1}^{\infty}a_{m}\sum_{i\in\sigma_{m}}\alpha_{i,m}\sum_{j\in E}Te_{i,j}\right)\right\Vert _{p}\\
 & = & \left|E\right|^{-1/p}\left\Vert T\right\Vert ^{-1}H_{p}^{-1}A_{p}\left(\int\left(\sum_{m=1}^{\infty}a_{m}^{2}u_{m}^{2}\left(E\right)\right)^{p/2}d\mu\right)^{1/p}\\
 & \geq & \left|E\right|^{-1/p}\left\Vert T\right\Vert ^{-1}H_{p}^{-1}A_{p}\left(\int\left(\sum_{m=1}^{\infty}a_{m}^{2}u_{m}^{2}\left(E\right)\chi_{B_{m}}\right)^{p/2}d\mu\right)^{1/p}\\
 & = & \left|E\right|^{-1/p}\left\Vert T\right\Vert ^{-1}H_{p}^{-1}A_{p}\left(\sum_{m=1}^{\infty}\left|a_{m}\right|^{p}\int_{B_{m}}u_{m}^{p}\left(E\right)d\mu\right)^{1/p}\\
 & \geq & \left|E\right|^{-1/p}\left\Vert T\right\Vert ^{-1}H_{p}^{-1}A_{p}\e_{0}\left(\sum_{m=1}^{\infty}\left|a_{m}\right|^{p}\right)^{1/p}.\end{array}\]
 \endpf

\begin{lm}
There are successive convex combinations $\nu_{k}\left(\cdot\right)$
of $\left\{ v_{i}^{2}\left(\cdot\right)\right\} $ such that for all
finite $E\subseteq\mathbb{N}$ $\nu_{k}\left(E\right)\underset{k\rightarrow\infty}{\rightarrow}\Lambda\left(E\right)$
in $L_{p/2}$. $\Lambda\left(E\right)$ is a $L_{p/2}^{+}$ additive
valued measure, $\Lambda$, on the finite subsets on $\mathbb{N}$.\end{lm}

\pf
{\bf Case 1} ($1<p<2$): Let $V=\left\{ \left(\sum_{n=1}^{\infty}\alpha_{n}^{2}v_{n}^{2}\right):\sum\alpha_{n}^{2}\leq1\right\} $.
Since $V$ is bounded in $L_{p/2}$, by a result of Nikishin \cite{NI}
for each $\varepsilon>0$ there is a set $D=D_{\varepsilon}\subset[0,1]$
of measure larger than $1-\varepsilon$ such that $\sup_{v\in V}\int_{D}vd\mu<\infty$.
As in the proof of \cite{JMST} Lemma 6.4 (or see Proposition 5.2 in \cite{DJS} for more details), we can find successive
convex combinations $\nu_{k}(\cdot)$ of the $v_{n}^{2}(\cdot)$ such
that $\nu_{k}(E){\bf 1}_{D}$ converges pointwise and in $L_{1}$
to $\Lambda(E){\bf 1}_{D}$ for every finite $E\subset\mathbb{N}$,
where $\Lambda\left(E\right)$ is $L_{0}^{+}$-valued (and $\Lambda{\bf 1}_{D}$
is $L_{1}^{+}$-valued). By passing to a subsequence of the $\nu_{k}$
and a simple diagonal argument we can find, for every $\varepsilon_{n}\to0$
a sequence of sets $D_{n}$ with $\mu(D_{n})>1-\varepsilon_{n}$ and
such that $\nu_{k}(E){\bf 1}_{D_{n}}$ converges, as $k\to\infty$,
pointwise and in $L_{1}$ to $\Lambda(E){\bf 1}_{D_{n}}$ for every
finite $E\subset\mathbb{N}$ and every $n$. In particular, $\nu_{k}(E)$
converges pointwise to $\Lambda(E)$ for every finite $E\subset\mathbb{N}$.
It remains to show that the convergence is also in $L_{p/2}$ (on
the whole interval). Since for each $E$, $\left\{ \nu_{k}\left(E\right)\right\} _{k\in\mathbb{N}}$
is $p/2$-equi-integrable, it follows that, given any $\delta>0$,
if $n$ is large enough $\int_{D_{n}^{c}}\nu_{k}(E)^{p/2}d\mu<\delta$
for all $k$. Consequently, also $\int_{D_{n}^{c}}\Lambda(E)^{p/2}d\mu\le\delta$
and

\[
\begin{array}{rcl}
\limsup_{k\to\infty}\int|\nu_{k}(E)-\Lambda(E)|^{p/2}d\mu & \leq & \limsup_{k\to\infty}\int_{D_{n}}|\nu_{k}(E)-\Lambda(E)|^{p/2}d\mu+2\delta\\
 & \leq & \limsup_{k\to\infty}\|(\nu_{k}(E)-\Lambda(E)){\bf 1}_{D_{n}}\|_{1}^{p/2}+2\delta\\
 & = & 2\delta.\end{array}\]

{\bf Case 2} ($2<p<\infty$):
\[
\left\Vert v_{i}^{2}\left(E\right)\right\Vert _{p/2}^{1/2}=\left\Vert v_{i}\left(E\right)\right\Vert _{p}\leq H_{p}\left\Vert \sum_{j\in E}Te_{i,j}\right\Vert _{p}\leq
H_{p}\left\Vert T\right\Vert \left\Vert \sum_{j\in E}e_{i,j}\right\Vert _{\ell_p(\ell_2)}=H_{p}\left\Vert T\right\Vert \left|E\right|^{1/p}.
\]
Thus, there exists a subsequence of $\left\{ v_{i}^{2}\right\} _{i=1}^{\infty}$
that converges weakly in $L_{p/2}\left[0,1\right]$ for every $E\in2^{\mathbb{N}}$.
Denote the limit by $\Lambda\left(E\right)$. By the reflexivity
of $L_{p/2}\left[0,1\right]$ and another diagonal argument, there
exists a sequence $\left\{ \sigma_{k}\right\} _{k=1}^{\infty}$ of
disjoint finite subsets of the integers and non-negative numbers $\left\{ \alpha_{i}\right\} _{i=1}^{\infty}$
such that $\sum_{i\in\sigma_{k}}\alpha_{i}^{2}=1$ and $\sum_{i\in\sigma_{k}}\alpha_{i}^{2}v_{i}^{2}\left(E\right)\rightarrow\Lambda\left(E\right)$
as $k\rightarrow\infty$ for all finite $E$ where the convergence
is in the $L_{p/2}\left[0,1\right]$ norm. Set $\nu_{k}\left(\cdot\right)=\sum_{i\in\sigma_{k}}\alpha_{i}^{2}v_{i}^{2}\left(\cdot\right)$.
$\nu_{k}$ is clearly additive on the finite subsets of $\mathbb{N}$.\endpf

It is clear that the sequence $\left\{ \Lambda\left(\left\{ j\right\} \right)\right\} _{j=1}^{\infty}$
is positively equivalent to the unit vector basis of $\ell_{p/2}$;
i.e., putting $\Lambda_{j}=\Lambda\left(\left\{ j\right\} \right)$,
$\left\Vert \sum a_{j}^{2}\Lambda_{j}\right\Vert _{p/2}\approx\left(\sum\left|a_{j}\right|^{p}\right)^{2/p}$
for all the sequences of coefficients $\left\{ a_{j}\right\} _{j=1}^{\infty}$.
Next we would like to improve the equivalence constant to be arbitrarily close
to 1 by blocking the $\Lambda_{j}$-s.

\begin{lm}\label{lm:almostdisjoint}
Let $1<p<\infty$, $\e>0$ and $\e_{k}\searrow0$. There
are successive disjoint $\sigma_{k}\subseteq\mathbb{N}$ $k=1,2,...$
and coefficients $\left\{ \alpha_{j}\right\} _{j=1}^{\infty}$ such
that putting $\phi_{k}=\sum_{j\in\sigma_{k}}\alpha_{j}^{2}\Lambda_{j}$,
for all $\left\{ a_{k}\right\} _{k=1}^{\infty}$,
\[
\left(\sum_{k=1}^{\infty}\left|a_{k}\right|^{p}\right)^{2/p}\leq\left\Vert \sum_{k=1}^{\infty}a_{k}^{2}\phi_{k}\right\Vert _{p/2}\leq\left(1+\e\right)\left(\sum_{k=1}^{\infty}\left|a_{k}\right|^{p}\right)^{2/p}.
\]
In addition, there are disjoint sets $\left\{ B_{k}\right\} _{k=1}^{\infty}$
such that $\left\Vert \phi_{k}{\bf 1}_{B_{k}^{c}}\right\Vert \leq\e_{k}$.
\end{lm}

\pf
By the subsequence splitting lemma (see e.g. \cite{AK}, Lemma 5.2.8) we may
assume, passing to a subsequence, that there are disjoint sets $\left\{ A_{j}\right\} _{j=1}^{\infty}$
such that $\left\{ \Lambda_{j}^{p/2}{\bf 1}_{A_{j}^{c}}\right\} _{j=1}^{\infty}$
is equi-integrable. For $p>2$ this already implies that $\int\Lambda_{j}^{p/2}{\bf 1}_{A_{j}^{c}}\rightarrow0$.
Otherwise, there is an $R>0$ and $\delta>0$ such that $\int\Lambda_{j}^{p/2}{\bf 1}_{\Lambda_{j}\leq R}\geq\delta$
for $j\in J$ for some infinite $J\subseteq\mathbb{N}$ and it would follow
that for $j\in J$ $\int\Lambda_{j}\geq\delta R^{1-p/2}$
and for all coefficients $\{a_i\}$, $\int\left(\sum a_{j}^{2}\Lambda_{j}\right)^{p/2}\geq\left(\int\sum a_{j}^{2}\Lambda_{j}\right)^{p/2}\geq\delta^{p/2}R^{\frac{p}{2}\left(1-\frac{p}{2}\right)}\left(\sum a_{j}^{2}\right)^{p/2}$,
contradicting the positive equivalence to the unit vector basis of $\ell_{p/2}$.
We can now take the $\sigma_{k}$-s to be singletons and $B_{k}$-s
to be some subsequence of the $A_{j}$-s.

For $1<p<2$, the proof is a bit more complicated. Set $h_{j}=\Lambda_{j}{\bf 1}_{A_{j}^{c}}$ and let $\delta_k\searrow0$ be a sequence to be determined later.
By the equi-integrability, there exists an $R>0$ such that
\[
\int\left|\frac{1}{n^{2/p}}\sum_{j=1}^{n}h_{j}{\bf 1}_{\left\{ h_{j}\geq R\right\} }\right|^{p/2}d\mu\leq\frac{1}{n}\sum_{i=1}^{n}\int h_{j}^{p/2}{\bf 1}_{\left\{ h_{j}\geq R\right\} }d\mu\leq\frac{\delta_{1}^{p/2}}{2}.
\]
We also have,
\[
\int\left|\frac{1}{n^{2/p}}\sum_{j=1}^{n}h_{j}{\bf 1}_{\left\{ h_{j}<R\right\} }\right|^{p/2}d\mu<\int\left|\frac{1}{n^{2/p}}\sum_{j=1}^{n}R\right|^{p/2}d\mu= R^{p/2}\frac{1}{n}n^{p/2}\leq\frac{\delta_{1}^{p/2}}{2}
\]
for $n$ sufficiently large. Set $m_{1}=\frac{1}{n^{2/p}}\sum_{j=1}^{n}\Lambda_{j}$
and $B_{1}=\cup_{j=1}^{n}A_{j}$. Then
\[
\begin{array}{rl}
\int\left|m_{1}{\bf 1}_{B_{1}^{c}}\right|^{p/2}d\mu&=
\int\left|\frac{1}{n^{2/p}}\sum_{j=1}^{n}\Lambda_{j}{\bf 1}_{B_{1}^{c}}\right|^{p/2}d\mu\\
&\leq
\int\left|\frac{1}{n^{2/p}}\sum_{j=1}^{n}\Lambda_{j}{\bf 1}_{A_{j}^{c}}\right|^{p/2}d\mu\leq\delta_{1}^{p/2}.
\end{array}
\]
Similarly, we find $m_{2}=\frac{1}{\left(n_{2}-n_{1}\right)^{2/p}}\sum_{j=n_{1}+1}^{n_{2}}\Lambda_{j}$
and $B_{2}=\cup_{j=n_{1}+1}^{n_{2}}A_{j}$ that satisfy
\[
\int\left|m_{2}{\bf 1}_{B_{2}^{c}}\right|^{p/2}d\mu\leq \delta_2^{p/2}.\]
 Continuing in this manner we find $m_i=\frac{1}{\left(n_{i}-n_{i-1}\right)^{2/p}}\sum_{j=n_{i-1}+1}^{n_{i}}\Lambda_{j}$ and disjoint $\{B_j\}$ satisfying $\int\left|m_{i}{\bf 1}_{B_{i}^{c}}\right|^{p/2}d\mu\leq \delta_i^{p/2}$. The sequence $\{\|m_k\|_{p/2}\}$ is bounded and bounded away from zero. Consequently, putting $\Phi_k=\frac{m_k}{\|m_k\|_{p/2}}$ we get that for an appropriate choice of the $\delta_k$-s, $\{\Phi_k,B_k\}$ satisfy the conclusion of the lemma. \endpf

{\bf Proof of Theorem \ref{thm:main}} In the following let $\sigma_k$ and $\alpha_j$ be as in the statement of Lemma \ref{lm:almostdisjoint}. Let $\e_{l,k}>0$ and for each $k$ let $\left\{ \sigma_{l,k}\right\} _{l=1}^{\infty}$
be successive finite subsets of $\mathbb{N}$ and $\left\{ \beta_{i,k}\right\} _{i\in\sigma_{l,k}}$
coefficients such that $\sum_{i\in\sigma_{l,k}}\beta_{i,k}^{2}=1$ and
$\left\Vert \sum_{i\in\sigma_{l,k}}\beta_{i,k}^{2}\nu_{i}^{2}\left(\sigma_{k}\right)-\Lambda\left(\sigma_{k}\right)\right\Vert _{p/2}<\e_{l,k}$.
Then
\[
\left\Vert \sum_{j\in\sigma_{k}}\alpha_{j}^{2}\sum_{i\in\sigma_{l,k}}\beta_{i,k}^{2}\nu_{i}^{2}\left(j\right)
-\sum_{j\in\sigma_{k}}\alpha_{j}^{2}\Lambda_{j}\right\Vert _{p/2} \leq
\left(\max_{j\in\sigma_{k}}\left|\alpha_{j}\right|^{p}\right)^{2/p}\e_{l,k}.
\]
Put $f_{l,k}=\sum_{j\in\sigma_{k}}\alpha_{j}\sum_{i\in\sigma_{l,k}}\beta_{i,k}e_{i,j}$.
Then, if the $\e_{l,k}$-s are small enough,
\[
\begin{array}{rcl}
\left\Vert \sum_{l,k}a_{l,k}f_{l,k}\right\Vert _{\ell_p(\ell_2)} & \approx & \left\Vert \sum_{l,k}a_{l,k}^{2}S^{2}\left(T\left(f_{l,k}\right)\right)\right\Vert _{p/2}^{1/2}\\
 & = & \left\Vert \sum_{l,k}a_{l,k}^{2}\sum_{j\in\sigma_{k}}\alpha_{j}^{2}\sum_{i\in\sigma_{l,k}}\beta_{i,k}^{2}\nu_{i}^{2}\left(j\right)\right\Vert _{p/2}^{1/2}\\
 & \overset{1+\e}{\approx} & \left\Vert \sum_{l,k}a_{l,k}^{2}\sum_{j\in\sigma_{k}}\alpha_{j}^{2}\Lambda_{j}\right\Vert _{p/2}^{1/2}\\
 & \overset{1+\e}{\approx} & \left(\sum_{k=1}^{\infty}\left(\sum_{l=1}^{\infty}a_{l,k}^{2}\right)^{p/2}\right)^{1/p}.\end{array}\]
This shows that $\left\{ f_{l,k}\right\} _{l=1,k=1}^{\infty,\infty}$
is equivalent to the unit vector basis of $\ell_p(\ell_2)$. The constant
of the equivalence depends however on $\left\Vert T\right\Vert ,\left\Vert T^{-1}\right\Vert $
and $p$. We next correct this: we can find, for each $k$, a block
basis $g_{u,k}=\sum_{l\in\tau_{u,k}}\gamma_{l,k}f_{l,k}$, $u=1,2,...$,
such that $\left\Vert \sum_{u=1}^{\infty}a_{u}g_{u,k}\right\Vert _{\ell_p(\ell_2)}\overset{1+\e}{\approx}\left(\sum a_{u}^{2}\right)^{1/2}$.\\
Since for each $k$ $g_{u,k}$ is supported only on columns in
$\sigma_{k}$, we get from that that $\left\{ g_{u,k}\right\} _{u,k}$
is $\left(1+\e\right)$-equivalent to the unit vector basis of $\ell_p(\ell_2)$.

Next we would like to apply a similar stabilization procedure to $\left\{ T\left(g_{u,k}\right)\right\} _{u,k}$.
This of course is more complicated since this sequence does not lie
in $\ell_p(\ell_2)$ any more.

Note that, by passing to further subsequences we may assume that
for some positive constant $b$ (depending on $T$), and for all $k$
and $\left\{ a_{n}\right\} $ with $\sum_{n=1}^{\infty}a_{n}^{2}=1$
\[
\left\Vert S^{2}\left(\sum_{n=1}^{\infty}a_{n}Tg_{n,k}\right)-b\phi_{k}\right\Vert _{p/2}<\e_{k}\]
 where $\e_{k}$ is the preassigned sequence appearing in Lemma \ref{lm:almostdisjoint}.
In particular,
\[
\begin{array}{rl}
&\left\Vert S^{2}\left(\sum_{n=1}^{\infty}a_{n}Tg_{n,k}\right){\bf 1}_{B_{k}^{c}}\right\Vert _{p/2}^{p/2}\\
&
\phantom{aaaaa}\leq
\left\Vert (S^{2}\left(\sum_{n=1}^{\infty}a_{n}Tg_{n,k}\right)-b\phi_{k}){\bf 1}_{B_{k}^{c}}\right\Vert _{p/2}^{p/2} + b^{p/2}\left\Vert \phi_{k}{\bf 1}_{B_{k}^{c}}\right\Vert _{p/2}^{p/2}\\
 & \phantom{aaaaa}\leq
\e_{k}^{p/2}+b^{p/2}\e_{k}^{p/2}.\end{array}
\]

This shows that any sequence of the form $\{S(\sum_{n=1}^\infty a_{n,k}Tg_{n,k})\}_{k=1}^\infty$ with, say, \\
$\sum_{n=1}^\infty a_{n,k}^2=1$ for all $k$ is essentially (with respect to the $L_p$ norm) disjointly supported. We would like to have a similar statement with the $S$ removed.
This is where the Haar functions (rather than some other unconditional basis of $L_p$) play a role. First we may assume that the sets
$B_{k}$ are each a finite union of dyadic intervals. Next notice that
for each $k$ if $l$ is large enough then $Tf_{l,k}$ is a dyadic
simple function such that each Haar function appearing in its expansion (with non zero coefficient) has support which is either
contained or disjoint of $B_{k}$. Consequently,
\begin{equation}\label{eq:comm}
S\left(Tf_{l,k}{\bf 1}_{B_{k}}\right)=S\left(Tf_{l,k}\right){\bf 1}_{B_{k}}\ \
\mbox{and}\ \ S\left(Tf_{l,k}{\bf 1}_{B_{k}^{c}}\right)=S\left(Tf_{l,k}\right){\bf 1}_{B_{k}^{c}}.
\end{equation}
In particular, this implies that for some constant $K_p$ (depending on $p$ only), for all $k$ and all $\left\{ a_{n}\right\} _{n=1}^{\infty}$,
\[
\begin{array}{rcl}
\left\Vert \left(\sum_{n=1}^{\infty}a_{n}Tg_{n,k}\right){\bf 1}_{B_{k}^{c}}\right\Vert _{p} & \leq & K_{p}\left\Vert S^{2}\left(\sum_{n=1}^{\infty}a_{n}Tg_{n,k}\right){\bf 1}_{B_{k}^c}\right\Vert _{p/2}^{1/2}\\
 & \leq & K_{p}\left(1+b^{p/2}\right)^{1/p}\e_{k}^{1/2}\left(\sum_{n=1}^{\infty}a_{n}^{2}\right)^{1/2}.
 \end{array}
 \]
We thus get that if $\e_{k}$ are chosen small enough, the
two sequences $\left\{ Tg_{n,k}{\bf 1}_{B_{k}}\right\} _{n,k=1}^{\infty}$ and $\left\{ Tg_{n,k}\right\} _{n,k=1}^{\infty}$
are small perturbations one of the other. Define $\tilde{T}:\ell_p(\ell_2)\rightarrow\left[\left(Tg_{n,k}\right){\bf 1}_{B_{k}}\right]_{n,k=1}^{\infty}$
by $\tilde{T}e_{n,k}=\left(Tg_{n,k}\right){\bf 1}_{B_{k}}$. The perturbation
above is such that if we show that for some subspace $X\subset \ell_p(\ell_2)$,
spanned by blocks $\left\{ u_{n,k}\right\} _{n,k=1}^{\infty}$ of $\left\{ e_{n,k}\right\} _{n,k=1}^{\infty}$
$1$-equivalent to the unit vector basis of $\ell_p(\ell_2)$, and such
that $\left\{ \tilde{T}u_{n,k}\right\} _{n,k=1}^{\infty}$ is \\
$\left(1+\e\right)$-equivalent
to a multiple of the unit vector basis  of $\ell_p(\ell_2)$ and $\left[\tilde{T}u_{n,k}\right]_{n,k=1}^{\infty}$
is 
$\left(1+\e\right)\gamma_{p/(p-1)}$-complemented in $L_{p}$, then the same
(with $1+2\e$ replacing $1+\e$) holds for $\left\{ Tu_{n,k}\right\} _{n,k=1}^{\infty}$.

The way we choose $u_{n,k}$ is similar to the way we have chosen
$g_{n,k}$: since for each $k$ $\left\{ Tg_{n,k}\right\} _{k=1}^{\infty}$
is equivalent to the unit vector basis  of $\ell_2$ we can find blocks $\left\{ u_{n,k}\right\} $
of $\left\{ g_{n,k}\right\} $ such that $\left\{ \left(Tu_{n,k}\right){\bf 1}_{B_{k}}\right\} $
(which are blocks of $\left\{ \left(Tg_{n,k}\right){\bf 1}_{B_{k}}\right\} _{n,k=1}^{\infty}$)
are $\left(1+\e\right)$ equivalent to the unit vector basis of $\ell_2$
and $\left(1+\e\right)\gamma_{p}$ complemented in $L_{p}$.
Call the projection $P_{k}$. Note also that we may assume that the
$\ell_2$ norm of the coefficients of $\left\{ Tu_{n,k}\right\} $ relative
to $\left\{ Tf_{n,k}\right\} $ all differ by a multiplicative constant
of at most $1+\e$. (This is important since we want $T_{\mid X}$
to be a $\left(1+\e\right)$-isomorphism). Finally, define $P:L_{p}\rightarrow\left[\left(Tu_{n,k}\right){\bf 1}_{B_{k}}\right]$
by $Pf=\sum P_{k}\left(f{\bf 1}_{B_{k}}\right)$.\endpf

\section{Concluding remarks}
 \noindent 1. We first remark that the constants $\gamma_p$ and $\gamma_{p/(p-1)}$ that appear in the statement of Theorem \ref{thm:main} (and also in the theorem in section \ref{sec:3}) are best possible. Actually, these are lower bounds on the norm of the best projection onto an (isometric) $\ell_2$ subspace of $L_p$. This was proved in $\cite{GLR}$.
 
 \noindent 2. If one wants to avoid the use of the material in section \ref{sec:3} at the price of getting worse constants, still depending only on $p$, one can use Theorem 3.1 in \cite{PR} instead.
 
 \noindent 3. As we already remarked, it would be nice if somebody comes up with a simpler proof of the theorem of section \ref{sec:3}. Maybe along the lines of \cite{Als}.
 
 \noindent 4. Since $\ell_2$ is not isomorphic to a complemented subspace of $L_1$ there is of course no complete analogue of Theorem \ref{thm:main} for $L_1$. However, the weaker statement that any isomorphism $T:\ell_1(\ell_2)  \to L_1$ stabilizes; i.e., is a $(1+\e)$-isomorphism when restricted to some subspace of $\ell_1(\ell_2)$  $(1+\e)$-isomorphic to $\ell_1(\ell_2)$ , is still possible.
Of course the proof as is written above makes an heavy use of the unconditionality of the Haar system and thus cannot be used. However, note that for most of the proof we could replace the Haar system with any unconditional basic sequence containing $T(\ell_1(\ell_2))$, even with constant depending on $T$. So we could use $T e_{i,j}$ as such a sequence. The problem is in (\ref{eq:comm}) where we use the ``eventual commutativity" of the square function operation and the restriction to a given dyadic set. This seems like a not very essential use of the Haar system but we couldn't overcome it easily and we leave it for future research.

\begin{tabular}{ll}
R. Levy&G. Schechtman\\
Department of Mathematics&Department of Mathematics\\
Weizmann Institute of Science&Weizmann Institute of Science\\
Rehovot, Israel&Rehovot, Israel\\
{\tt  ran.levy@weizmann.ac.il }
&{\tt gideon@weizmann.ac.il}\\
&\\
R. Levy current address:&\\
IBM Research Lab&\\
Haifa, Israel&\\
{\tt ranl@il.ibm.com}
\end{tabular}


\begin{thebibliography}{99}
\bibitem[Ald]{Ald} D. J. Aldous, Subspaces of $L^{1}$, via Random
Measures, Trans. of Amer. Math. Soc., 267, No. 2 (Oct., 1981), pp.
445-463.

\bibitem[Als]{Als} D. E. Alspach, Good $\ell_2$-subspaces of $L_{p}$,
$p>2$, Banach J. Math. Anal. 3 (2009), no. 2, 49\textendash{}54

\bibitem[AK]{AK} F. Albiac and N.J. Kalton, Topics in Banach space
theory, Graduate Texts in Mathematics, 233, Springer, New York, 2006.

%\bibitem[B]{B} D. L. Burkholder, Martingale Transforms, Ann. of Prob.,
%9, No. 6 (Dec., 1981), pp. 997-1011.

\bibitem[DJS]{DJS} D. Dosev, W.B. Johnson and G. Schechtman, Comutators on $L_p$, $1\le p<\infty$,
 http://arXiv.org/abs/1102.0137

\bibitem[ES]{ES}P. Enflo and T. Starbird, subspaces of $L_{1}$ containing
$L_{1}$, Studia Math. 65(1979), 203-225.

\bibitem[GLR]{GLR} Y. Gordon, D. R. Lewis and J. R. Retherford, Banach ideals of operators with applications, J. Functional Analysis 14 (1973), 85–129.

\bibitem[HOS]{HOS} R. Haydon, E. Odell and T. Schlumprecht, Small
Subspaces of $L_{p}$, Annals of Math., to appear.

\bibitem[JMST]{JMST} W.B. Johnson, B. Maurey, G. Shechtman and L.
Tzafriri, Symmetric structures in Banach spaces, Memoirs Amer. Math.
Soc., 217 (1979).

\bibitem[KM]{KM} J.L. Krivine and B. Maurey, Espaces de Banach stables, (French), Israel J. Math. 39 (1981), no. 4, 273-295.

\bibitem[KP]{KP} M. I. Kadec and A. Pelczynski, Bases, lacunary sequences
and complemented subspaces in the spaces $L_{p}$, Studia Math. 21
(1962), pp. 161-176.

\bibitem[LT-I]{LT-I} J. Lindenstrauss and L. Tzafriri, Classical Banach spaces I, Sequence spaces, Springer-Verlag, Berlin-New York, 1977.

\bibitem[LT-II]{LT-II} J. Lindenstrauss and L. Tzafriri, Classical
Banach spaces II, Function spaces, Springer-Verlag, New York, 1979.

\bibitem[NI]{NI} E.M. Nikishin, Resonance theorems and superlinear
operators, Russ. Math. Surv. 25, no. 6, (1970), pp. 125-187.

\bibitem[PR]{PR} A. Pelczynski and H. P. Rosenthal, Localization
techniques in $L_{p}$ spaces, Studia Math. 52, (1975), pp. 263\textendash{}289.

\bibitem[R]{R}H. P. Rosenthal, Embeddings of $L_{1}$ in $L_{1}$,
Contemporary Math. 29(1984), 335-349.
\end{thebibliography}
\end{document}